\documentclass[10pt]{article}
\usepackage[utf8]{inputenc}

\usepackage{amsfonts}
\usepackage{amsmath}
\usepackage{amssymb}
\usepackage{afterpage}
\usepackage{amsthm}
\usepackage{mathtools}
\usepackage{esint}
\usepackage{mathrsfs}
\usepackage{bbm}
\usepackage{comment, todonotes}
\usepackage{protosem}
\usepackage{longtable}

\usepackage{mathabx}
\usepackage{graphicx}
\usepackage[font=footnotesize, labelfont=bf]{caption}
\usepackage[font=footnotesize, labelfont=bf]{subcaption}

\definecolor{labelkey}{rgb}{0,0,1}

\usepackage{fullpage}
\usepackage{enumitem}
\usepackage{hyperref,xcolor,fullpage}

\usepackage{todonotes}

\allowdisplaybreaks % allows break in multiline formulas

\renewcommand*{\div}{\ensuremath{\mathrm{div}}}
\newcommand{\eps}{\varepsilon}

\newcommand{\p}{\partial}

\def\les{\lesssim}
\newcommand{\eye}{\textproto{o}}

\newcommand{\abs}[1]{\left|#1\right|}

\DeclarePairedDelimiterX{\ip}[2]{\langle}{\rangle}{#1, #2}

\numberwithin{equation}{section}

\title{Smooth self-similar imploding profiles to 3D compressible Euler\\
{\normalfont \normalsize \sffamily  \textcolor{black}{This review article is dedicated to Constantine Dafermos’ 80th birthday}}}
\date{}
\author{Tristan Buckmaster\thanks{\footnotesize Department of Mathematics, University of Maryland, College Park, MD 
 \href{tristanb@umd.edu}{{tristanb@umd.edu}}} , Gonzalo Cao-Labora\thanks{\footnotesize Department of Mathematics, 
Massachusetts Institute of Technology, Cambridge, MA \href{gcaol@mit.edu}{gcaol@mit.edu}}\, and 
Javier G\'omez-Serrano\thanks{\footnotesize Department of Mathematics,
Brown University,
Providence RI %\&         	
 \href{javier\_gomez\_serrano@brown.edu}{javier\_gomez\_serrano@brown.edu} %/ \href{jgomezserrano@ub.edu}{jgomezserrano@ub.edu}
 }
}

\begin{document}
\maketitle

\begin{abstract}
The aim of this note is to present the recent results in \cite{implosion}, concerning the existence of ``imploding singularities'' for the 3D isentropic compressible Euler and Navier-Stokes equations. Our work builds upon the pioneering work of Merle, Rapha\"el, Rodnianski and Szeftel \cite{MeRaRoSz19a,MeRaRoSz19c,MeRaRoSz19b} and proves the existence of self-similar profiles for  \emph{all} adiabatic exponents $\gamma>1$ in the case of Euler; as well as proving asymptotic self-similar blow-up for $\gamma=\frac75$ in the case of Navier-Stokes. Importantly, for the Navier-Stokes equation, the  solution is  constructed to have density bounded away from zero and constant at infinity,  the first example of blow-up in such a setting. For simplicity, we will focus our exposition on the compressible Euler equations.
\end{abstract}

\setcounter{tocdepth}{1}

%%% INTRODUCTION %%%%

\section{Introduction} \label{sec:intro}
The compressible Euler equations  describe the conservation of mass, momentum, and energy in a fluid, and are important in many fields, including aerodynamics and astrophysics. In this review, we present recent developments regarding the existence of smooth imploding solutions for the compressible Euler equations.
The full compressible Euler equations take the form
\begin{alignat*}{2}
{ \p_t( \rho ~u) +  \div (\rho u\otimes u+  p~{\rm Id}) \, }&{ =0 }\\
{  \p_t \rho + \operatorname{div} (\rho u) \, } &{  =0} \\
{ \partial_t E  +  \div((p+ E)u) } & { =0} 
\end{alignat*}
where $u$ is the velocity, $\rho$ is the density, $p$ is the pressure, and $E$ is the energy. The equations describe the  conservation of momentum, mass, and energy in a fluid, respectively.
The pressure is given by  the ideal gas law
\[p = (\gamma-1)(E-\frac12\rho\abs{u}^2 )=\frac1\gamma\rho^\gamma e^S\,,\]
 for the adiabatic exponent $\gamma>1$. The sound speed is given by $c=\sqrt{\frac{\gamma p}{\rho}}$. 

\subsection{Shock waves}

Before discussing implosion in detail, let us first describe the classical problem of shock waves, which can be seen as a prototypical singularity in the context of Euler's equations.
A shock wave occurs when the speed of a disturbance exceeds the local speed of sound. A fundamental problem in the the mathematical theory of compressible fluids is to provide a complete description of shock formation and development. In particular, one is interested in a complete description of the evolution of a smooth solution up until the point of singularity, and the shock, a co-dimension 1 space-time hypersurface, that proceeds the initial singularity.

The earliest rigorous result regarding shock wave formation traces back to the work of Lax \cite{Lax1964} in the 1D setting.  Generalizations and improvements of Lax’s result were obtained by John \cite{John1974} and Liu \cite{Li1979}, for the 1D Euler equations. See the book of Dafermos \cite{Da2010} for a more extensive
bibliography of 1D results.

Sideris \cite{Si1997} proved the existence of finite time singularities in 2D and 3D. The result of Sideris proves that a singularity occurs; however, it is not ascertained what form such a singularity takes.  Christodoulou \cite{Ch2007} and Christodoulou-Miao \cite{Ch2019} demonstrated the formation of shocks for 3D isentropic, irrotational fluids in the relativistic and non-relativistic settings respectively. Luk and Speck built on this work to handle the case of shock formation for 2D isentropic fluids with non-trivial vorticity \cite{LuSp2018}. The first author, together with Shkoller, and Vicol, employing a different approach, resolved the shock formation problem in the most general setting of full 3D compressible Euler \cite{buckmaster2020formation,buckmaster2020shock} (cf. \cite{2021arXiv210703426L}). The work \cite{buckmaster2020formation,buckmaster2020shock}, together with the prior work \cite{BuShVi2019} of the same authors were the first to isolate the self-similar profile of the initial singularity that precedes the development of shock waves. In particular, the works demonstrated that the asymptotic self-similar profile of the singularity is described by self-similar solutions to the Burgers' equation. 
More recently, Abbrescia and Speck \cite{2022arXiv220707107A} and Shkoller and Vicol \cite{ShVi22} have studied the problem of maximal development of shock waves.

With regards to shock development in one spatial dimension, global unique weak solutions satisfying the Rankine-Hugoniot conditions have been established (see \cite{glimm65,MR330788,MR1188562,Bressan1995}), but these methods neither provide a precise description of the shock front nor detect weak discontinuities: characteristic surfaces conjectured by Landau and Lifshitz \cite{landau_fluid_1987}. In multiple dimensions, Majda \cite{Ma1983b} studied the short-time evolution of the shock front starting from discontinuous initial data, which is smooth on either side of the shock front. 
This framework does not cover the shock development problem, where one must evolve from H\"older continuous pre-shock data and weak discontinuities may form. 
For the one-dimensional $p$-system (which models 1D isentropic Euler), Lebaud \cite{MR1309163} was the first to prove shock formation and development in her thesis work (cf. \cite{MR1860832,MR1897395}). In the case of the non-isentropic $3 \times 3$ Euler equations in spherical symmetry, shock formation and development were first established by Yin \cite{MR2085314}. Independently, Christodoulou and Lisibach \cite{MR3489205} demonstrated shock development for the barotropic Euler equations in spherical symmetry. The use of the isentropic model or the assumption of irrotational flow in higher dimensions has been referred to as \emph{restricted shock development} because it cannot produce weak solutions to the Euler equations. Christodoulou \cite{Ch2019} has also established restricted shock development for the irrotational and isentropic Euler equations in three dimensions outside of symmetry.

\begin{figure}[htp]
\centering
 \includegraphics[width=0.9\linewidth]{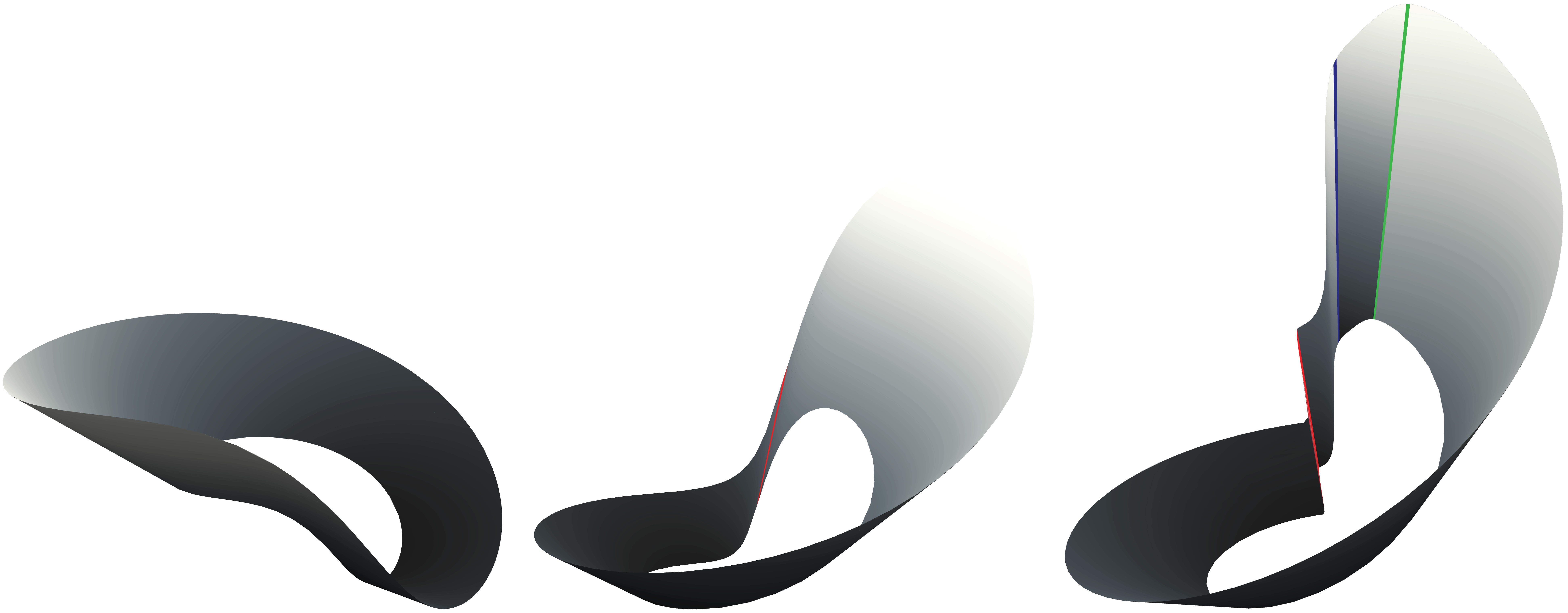}
\caption{\footnotesize The image represents the density restricted to the annular region $r\in[1,2]$. The first image is of the initial data, the second is the preshock, and the third is of the developed shock. The red, green and blue curves represent the shock curve,  weak rarefaction wave and  weak contact discontinuity respectively.}
\label{fig:headache}
\end{figure}

In \cite{BuDrShVi2021}, the first author, Drivas, Shkoller and Vicol consider the shock development problem for 2D compressible Euler under azimuthal symmetry (see \cite{ICM} for a recent review article). The work provides the first full description of shock development; in particular, in addition to describing the shock front, \cite{BuDrShVi2021} gives the first detailed description of the weak discontinuities of Landau and Lifshitz \cite{landau_fluid_1987} (see Figure \ref{fig:headache}).

\subsection{Implosion}

While shock waves are a common and potentially the only stable form of singularity for the Euler equations, other types of singularities can arise from smooth initial data. It is a fundamentally interesting problem to classify these forms of singularities, both from a mathematics and physics perspective.

Guderley's classical work \cite{MR8522} (cf.\ \cite{chisnell_1998,MeyerterVehn1982}) constructed the first examples of non-smooth imploding solutions. Very recently,  Merle, Rapha\"el, Rodnianski, and Szeftel rigorously proved the existence of \textit{smooth} radially symmetric imploding solutions to the isentropic compressible Euler equations \cite{MeRaRoSz19a}:
\begin{equation}\label{eq:Euler}
 \begin{split} 
\partial_t (\rho u) + \div (\rho u \otimes u) + \nabla p(\rho) &= 0 \,, \\
\partial_t \rho  +  \div (\rho u)&=0 \,, 
\end{split}
\end{equation}
where here $p(\rho) = \frac{1}{\gamma}\rho^{\gamma}$ for $\gamma > 1$.
Specifically, for almost every $\gamma>1$, they showed the existence of a countably infinite sequence of self-similar solutions to \eqref{eq:Euler}. These solutions exhibit blow-up of both the velocity and density at the origin. The condition on $\gamma$ is related to the non-vanishing of an analytic function. The case $\gamma=5/3$, which describes monatomic gases, is specifically ruled out.

The form of the singularity discovered in \cite{MeRaRoSz19a} is fundamentally new. The authors also used these solutions to prove finite-time blow-up for the defocusing, supercritical, nonlinear Schrödinger equation \cite{MeRaRoSz19b}, solving a significant open problem in the field. Additionally, the solutions were used as  a basis to construct asymptotically self-similar solutions to the three-dimensional isentropic compressible Navier-Stokes equations with density-independent viscosity \cite{MeRaRoSz19c}, given by
\begin{equation}\label{eq:NS}
 \begin{split} 
\partial_t (\rho u) + \div (\rho u \otimes u) + \nabla p(\rho)-\mu_1\Delta u-(\mu_1+\mu_2)\nabla \div u &= 0 \,, \\
\partial_t \rho  +  \div (\rho u)&=0 \,, 
\end{split}
\end{equation}
where $(\mu_1,\mu_2)$ are the Lamé viscosity coefficients, with $\mu_1>0$ and $2\mu_1+\mu_2> 0$. Prior to this result, Xin \cite{Xin} showed the existence of blow-up solutions for initial data with compact density, and  Rozanova \cite{ROZANOVA20081762} demonstrated the existence of blow-up solutions for rapidly (polinomially) decaying density. Unlike \cite{MeRaRoSz19c}, neither \cite{Xin} nor \cite{ROZANOVA20081762} provide a description of the singularity that occurs. The result \cite{MeRaRoSz19c} further weakens the decay required on the density leading to singularity formation. To rule out the role of vacuum at spatial infinity in the singularity formation, one however would prefer such solutions to be constructed from initial data that has non-vanishing, constant density at infinity. See also the recent numerical work by Biasi \cite{Biasi21}.

The papers \cite{MeRaRoSz19a} and \cite{MeRaRoSz19c} left open two fundamental questions:
\begin{enumerate}
\item
Do imploding solutions for the Euler equations exist for any value of $\gamma$ greater than 1?
\item Is it possible to create imploding solutions to the Navier-Stokes equation with an initial density that is constant at infinity?
\end{enumerate}

In \cite{implosion}, we resolved both of these questions. We showed that for all  $\gamma>1$ there exist self-similar imploding solutions. For the case of diatomic gases, $\gamma=\frac75$, we showed there exists an infinite sequence of self-similar imploding solutions. The paper \cite{implosion} also provides simplified proofs of  linear stability and non-linear stability, leading to the proof of asymptotically self-similar imploding solutions to the Navier-Stokes equations for $\gamma=\frac75$. The initial data for such solutions are chosen to have constant non-zero density at infinity -- the first example of such initial data leading to blow-up for the Navier-Stokes equations. The focus of this article will be on the former result.

\section{Reduction to an autonomous ODE}

\begin{figure}
\centering
  \includegraphics[width=.5\linewidth]{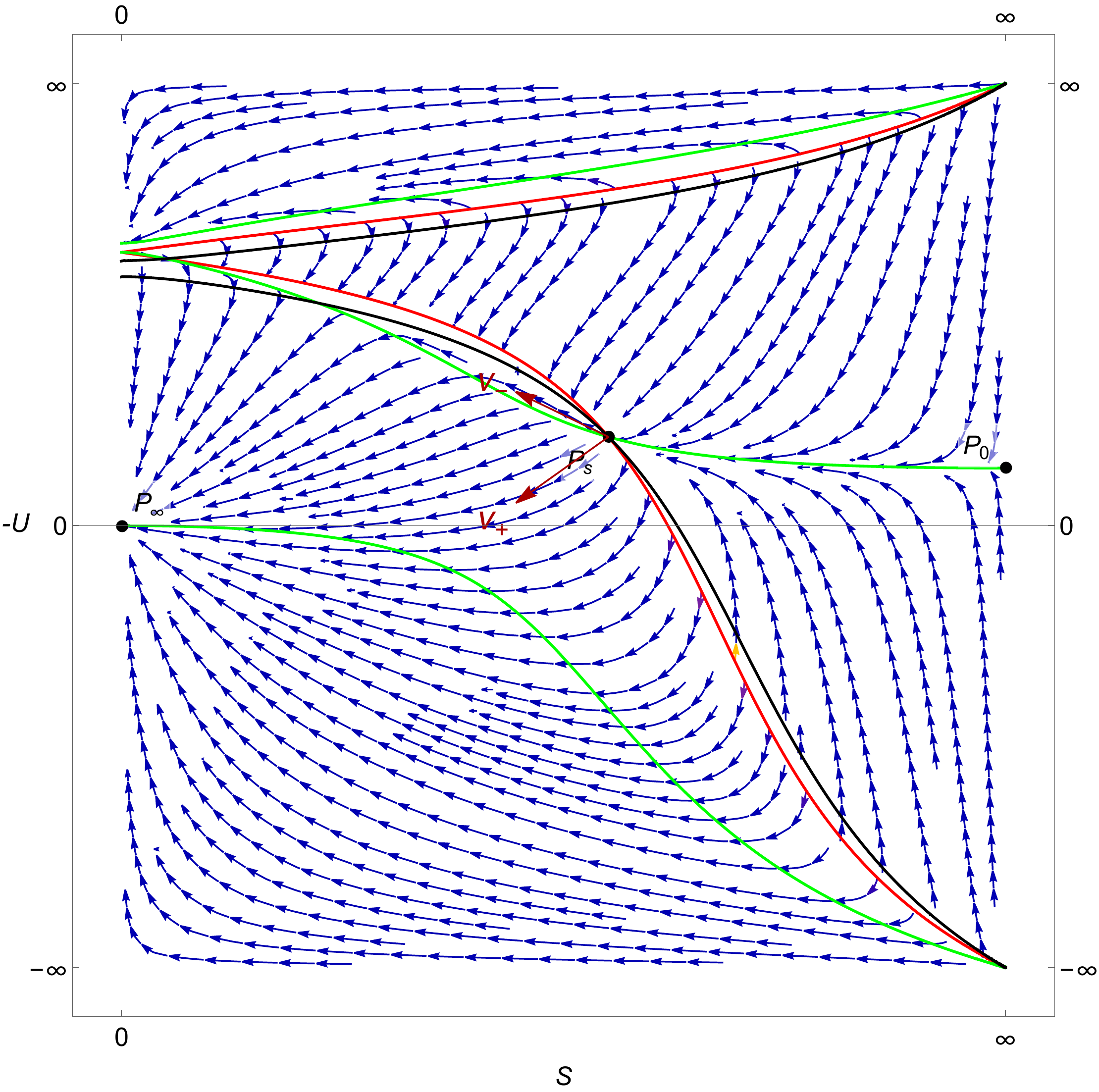}
  \captionof{figure}{\footnotesize Imploding solutions in $(U,S)$ variables. Note that a singular coordinate change has been made in order to compactify the $(U,S)$ coordinates.}
  \label{fig:US}
\end{figure}

Let us rewrite \eqref{eq:Euler} in radial form:
\begin{equation}\label{eq:wombat}
 \p_t u + u \p_R u + \frac{1}{\gamma\rho}\p_R \rho^{\gamma}=0\quad\mbox{and}\quad
 \p_t \rho + \frac{1}{R^2}\p_R(R^2 \rho u)=0\,,
\end{equation}
 where for matters of simplicity, we restricted the problem to three dimensions. Letting $\alpha = \frac{\gamma-1}{2}$, we define the rescaled sound speed: $\sigma=\frac{1}{\alpha}\rho^{\alpha}$. Then, we make the following self-similar anzatz
\[
  u(R,t)=r^{-1}\tfrac{R}{T-t} U(\log(\tfrac{R}{(T-t)^{\frac1r}}))\quad\mbox{and}\quad
 \sigma(R,t)=\alpha^{-\tfrac12}r^{-1} \tfrac{R}{T-t} S(\log(\tfrac{R}{(T-t)^{\frac1r}}))\,,
\]
 where here $r$ is a self-similar scaling parameter to be determined. Defining the self-similar variable $\xi=\log(\tfrac{R}{(T-t)^{\frac1r}})$, then \eqref{eq:wombat} reduces to an autonomous system of the form
 \begin{align}\label{eq:DS}
\tfrac{dU}{d\xi}  = \tfrac{N_U (U, S)}{D (U, S)},\quad\mbox{and}\quad
\tfrac{dS}{d\xi}  = \tfrac{N_S (U, S)}{D (U, S)}\,.
\end{align}
For, $\gamma=\frac75$ and $r=1.079404$, the phase portrait is shown in Figure \ref{fig:US}, where $D$, $N_U$, and $N_S$ are represented by red, green, and black curves, respectively. The point labeled $P_s$ is a special type of singular point for the dynamic system described in equation \eqref{eq:DS}. There are two smooth integral curves that pass through $P_s$, one tangent to the direction $\nu_-$ and the other tangent to $\nu_+$. The curve that is tangent to $\nu_+$ corresponds to the Guderley solution, while the curve tangent to $\nu_-$ corresponds to the solution found in \cite{MeRaRoSz19a}. To create a globally defined self-similar solution, we need to find an integral curve that connects the points $P_0$ and $P_\infty$ through $P_s$. It is not possible to do this using a continuous integral curve with the Guderley solution. However, by introducing a shock discontinuity, we can jump from one point in the phase portrait to another and create a globally defined self-similar solution. In \cite{MeRaRoSz19a}, by means of choosing distinguished values of the self-similar scaling parameter $r$, the authors overcame the challenge that the smooth integral curve tangent to $\nu_-$ generally does not connect $P_0$ to $P_\infty$, but rather intersects the sonic line $D=0$ at a point other than $P_s$, resulting in a solution that is not globally defined. 

Motivated by the works \cite{BuShVi2019,buckmaster2020formation,buckmaster2020shock}, it is helpful to rewrite the system in terms of its Riemann invariants
 \begin{equation}\label{eq:Riemann:invariants}
 w=u+\sigma\quad\mbox{and}\quad z=u-\sigma \,
 \end{equation}
so that
\[u=\frac{w+z}{2} \quad\mbox{and}\quad\sigma=\frac{w-z}{2}\,.\]
One can now diagonalize \eqref{eq:wombat} in terms of $w$ and $z$, in order to rewrite \eqref{eq:wombat}  as a nonlinear transport equation
 \begin{align}\label{eq:Euler:Riemann}\begin{split}
 \p_t w + \frac12(w+z+\alpha(w-z)) \p_R w + \frac{\alpha}{2R}(w^2-z^2)&=0\\
  \p_t z + \frac12(w+z-\alpha(w-z)) \p_R z - \frac{\alpha}{2R}(w^2-z^2)&=0\,.
  \end{split}
 \end{align}
Employing the self-similar ansatz
  \begin{align}\label{eq:ansatz}
  \begin{split}
 w(R,t)=\frac{1}{r} \cdot\frac{R}{T-t} W( \xi)\quad\mbox{and}\quad
 z(R,t)=\frac{1}{r} \cdot \frac{R}{T-t} Z( \xi)
 \end{split}
 \end{align}
 where we recall $\xi=\log(\tfrac{R}{(T-t)^{\frac1r}})$, then we obtain \begin{align} \begin{split} \label{eq:mainother}
(r+\frac12((1+2\alpha)W+(1-\alpha)Z))W+(1+\frac12(W+Z+\alpha(W-Z)))\partial_{\xi}  W  - \frac{\alpha}{2}Z^2&=0\\
(r+\frac12((1-\alpha)W+(1+2\alpha)Z))Z+(1+\frac12(W+Z-\alpha(W-Z)))\partial_{\xi}  Z  - \frac{\alpha}{2}W^2&=0\,
 \end{split} \end{align}
 
Rearranging, we obtain the autonomous system
\begin{equation}\label{eq:wz:ODE}
\begin{split}
\partial_{ \xi }  W&= \frac{-(r+\frac12((1+2\alpha)W+(1-\alpha)Z))W+ \frac{\alpha}{2}Z^2}{1+\frac12(W+Z+\alpha(W-Z))}=\frac{N_W}{D_W},\\
\partial_{ \xi } Z&=\frac{-(r+\frac12((1-\alpha)W+(1+2\alpha)Z))Z+\frac{\alpha}{2}W^2}{1+\frac12(W+Z-\alpha(W-Z))}=\frac{N_Z}{D_Z}.
\end{split}
\end{equation}

\begin{figure}
\centering
  \includegraphics[width=.5\linewidth]{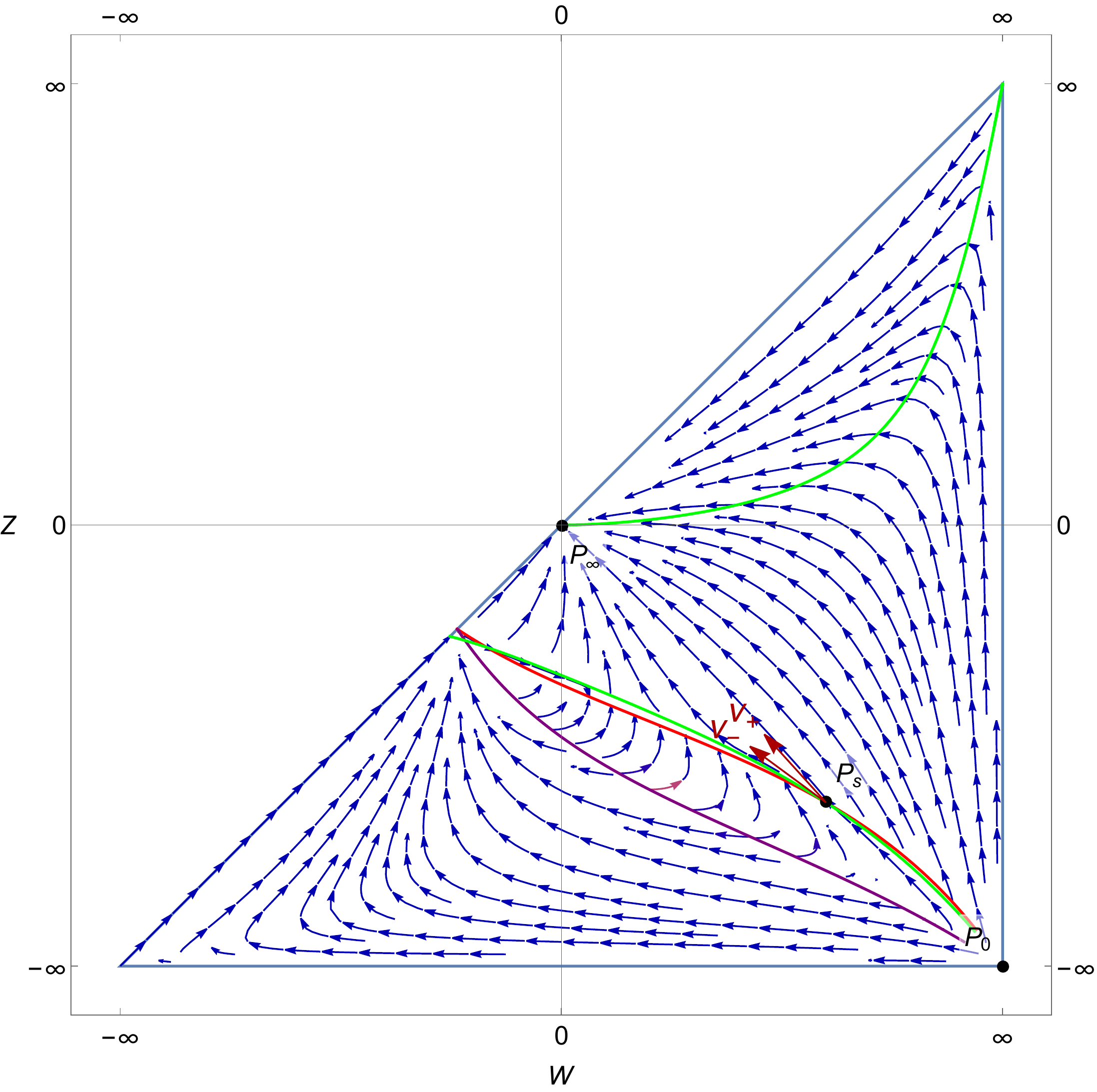}
  \captionof{figure}{\footnotesize Imploding solutions in $(W,Z)$ variables. Note that a singular coordinate change has been made in order to compactify the $(W,Z)$ coordinates.}
  \label{fig:WZ}
\end{figure}

In Figure \ref{fig:WZ}, the phase portrait for the region where the density is positive ($W-Z>0$) is shown. The red, purple, and green lines represent $D_Z=0$, $D_W=0$, and $N_Z=0$, respectively. One key difference between this system and \eqref{eq:DS} is that the denominator $D_W$ does not vanish at the point $P_s$, which simplifies the analysis in the area around $P_s$. The variables $(W,Z)$ provide a geometric understanding of the imploding solution in terms of the trajectories of the $W$ and $Z$ waves. $P_s$ is an unstable fixed point for the trajectories of $Z$-waves and divides space into an interior region (the backward acoustic cone emanating from the singular point) and an exterior region. $Z$-waves in the exterior region cannot enter the interior region, while $Z$-waves in the interior region cross the origin to become $W$-waves, then cross $P_s$ and travel to the exterior region. Since the system in \eqref{eq:wz:ODE} is autonomous, we can choose the location $\xi = 0$ to be where the solution crosses $P_s$.

The key steps to constructing a smooth integral curve from $P_0$ to $P_\infty$ are:
\begin{enumerate}
\item Apply a careful local analysis of the behavior of the smooth solution tangent to $\nu_-$. In particular, show that the solutions \emph{wiggle} in a certain manner with a continuous change in the self-similar parameter $r$.
\item Demonstrate that such a wiggling phenomenon combined with barrier arguments and continuity leads to a smooth solution connecting $P_0$ and $P_\infty$.
\end{enumerate}

 \section{Local analysis of $P_s$}
 
 To better understand $P_s$, we can recast the system using the variable $\xi\mapsto \psi$, where $\p_\psi=D_WD_Z\p_\xi$. For simplicity, let us focus on the case where $\gamma=\frac53$.
We obtain the ODE:
\begin{equation}\label{eq:wz:ODE2} 
\partial_{ \psi}  W=-\frac{1}{18} (3 + W + 2 Z) (6 r W + 5 W^2 + 2 W Z - Z^2)\quad\mbox{and}\quad
\partial_{ \psi } Z=\frac{1}{18} (3 + 2 W + Z) (W^2 - 2 (3 r + W) Z - 5 Z^2)\,,
\end{equation}
for which $P_s$ is a stable stationary point. Let $0<\lambda_-<\lambda_+$ be the eigenvalues of the resulting system of the Jacobian matrix at $P_s$. We let $k$ denote the ratio of the two eigenvalues:
\begin{equation}\label{eq:k:def}
k=\frac{\lambda_-}{\lambda_+}=\frac{r-2 - \sqrt{2r-2}}{r-2 + \sqrt{2r-2}}\,.
\end{equation}
The directions $\nu_-$ and $\nu_+$ defined earlier (the directions of the two smooth integral curves passing through $P_s$) are also the eigenvectors of the Jacobian of \eqref{eq:wz:ODE2} that correspond to the eigenvalues $\lambda_-$ and $\lambda_+$ respectively. We will focus on the smooth solutions of \eqref{eq:wz:ODE} that have tangents parallel to $\nu_-$. These two directions are shown in Figure \ref{fig:WZ}.

In the range $1 < r < r^\ast (= 3 - \sqrt{3}$ for $\gamma = \frac{5}{3}$), $k$ is a monotonically increasing function of $r$ that approaches infinity as $r$ approaches $r^\ast$. The smooth solution passing through point $P_s$ can be expressed as a Taylor series around $\xi = 0$ in the form $(W(\xi), Z(\xi)) = \sum_{n=0}^\infty \frac{\xi^n}{n!}(W_n, Z_n)$. The Taylor coefficients of $D_\circ(W, Z)$ and $N_\circ(W, Z)$ are denoted by $D_{\circ, n}$ and $N_{\circ, n}$, respectively. For $n \geq 2$, the following equations hold:
\begin{align} \begin{split} \label{eq:recurrence}
 D_{W, 0} W_n &= N_{W, n-1} - \sum_{j=0}^{n-2} \binom{n-1}{j} D_{W, n-1-j} W_{j+1},\\
Z_n D_{Z, 1} (n-k) &= - \sum_{j=1}^{n-2} \binom{n}{j} D_{Z, n-j} Z_{j+1} 
+\left( N_{Z, n} - (\p_Z N_Z (P_2) ) Z_n\right) + Z_1 \left( - D_{Z, n} + Z_n \p_Z D_Z (P_2) \right).
\end{split} \end{align}
By choosing $(W_1, Z_1)$ to align with $\nu_-$, these equations can be used to iteratively solve for a power series that describes the smooth solution tangent to $\nu_-$ at $P_s$ in a small neighborhood of $P_s$. Note that the right-hand side of the second equation does not depend on $Z_n$.

For any positive integer $j$, we define $r_j$ such that $j = k(r_j)$. It can be observed the expression for $Z_n$ in \eqref{eq:recurrence} becomes singular as $k(r)$ approaches $n$ and changes sign at $k(r) = n$. This causes the integral curve of the smooth solution to exhibit a \emph{wiggling} effect, which allows us to show\footnote{We believe this is true for all $\gamma > 1$ and every odd $n$.} that for $\gamma > 1$ and $n = 3$:
\begin{enumerate}
\item\label{pt:1} For $r \in (r_n, r_{n+1})$, the solution to the left of $P_s$ approaches $P_\infty$ as $\xi$ goes to infinity.
\item\label{pt:2} For $r = r_n + \eps$, the solution to the right of $P_s$ intersects the line $D_W = 0$.
\item\label{pt:3} For $r = r_{n+1} - \eps$, the solution to the right of $P_s$ intersects the line $D_Z = 0$.
\end{enumerate}
If we can demonstrate points \ref{pt:2} and  \ref{pt:3}, then using a simple shooting argument, we can conclude that there exists an $r$ within the range $(r_n, r_{n+1})$ such that the solution curve connects $P_s$ to $P_0$. Additionally, point \ref{pt:1} implies that the solution curve also connects $P_s$ to $P_\infty$.
 
\begin{figure}
\centering
  \includegraphics[width=.45\linewidth]{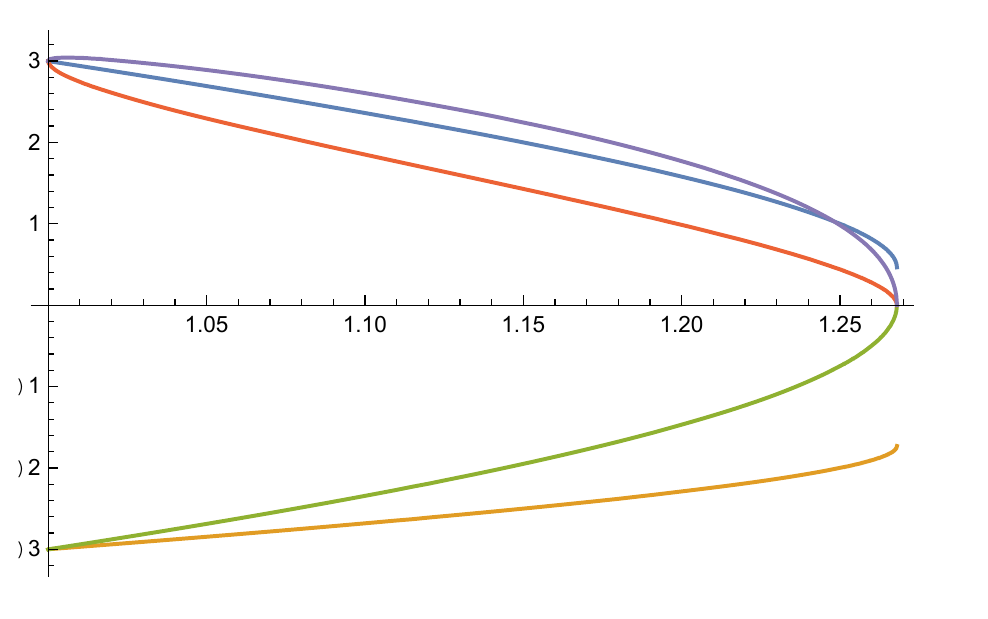}
  \includegraphics[width=.45\linewidth]{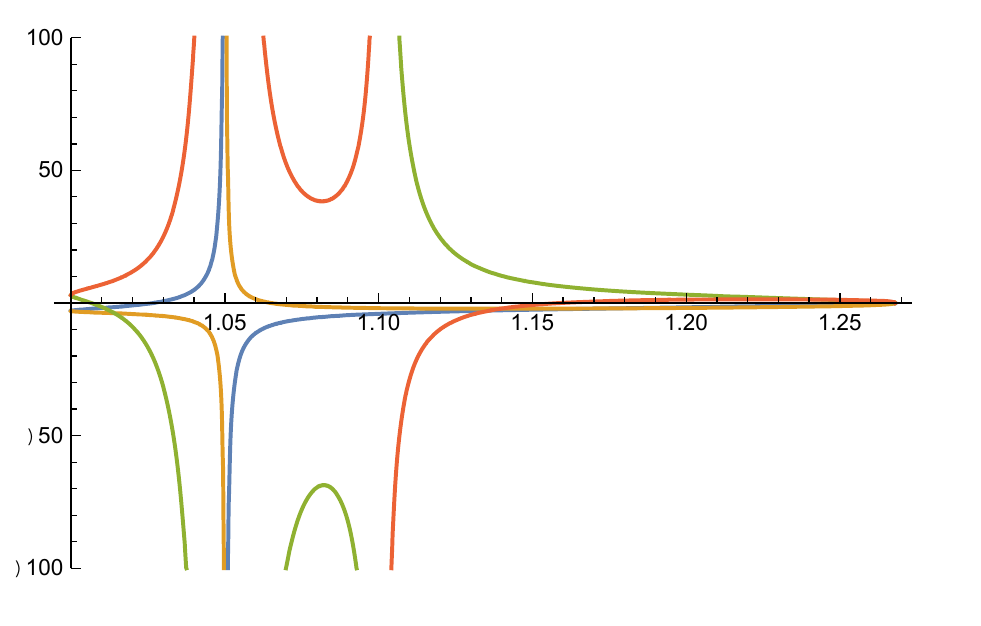}
  \includegraphics[width=.45\linewidth]{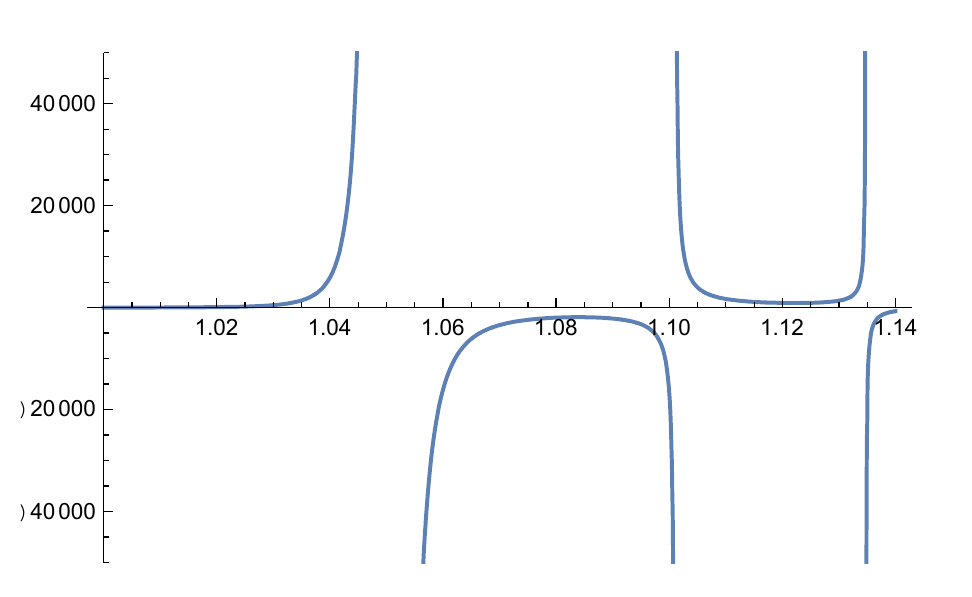}
  \captionof{figure}{\footnotesize In the first plot, we show $W_0$ (blue), $Z_0$ (orange), $W_1$ (green), $Z_1$ (red) and $W_2$ magenta. In the second plot, we show $Z_2$ (blue), $W_3$ (orange), $Z_3$ (green) and $W_4$ (red). In the last plot, we show $Z_4$. The singularities happen at $r = r_2 = 11-3\sqrt{11} \approx 1.05013$, $r = r_3 = 6-2\sqrt{6} \approx 1.10102$ and $r = r_4 = \frac19 (43-5\sqrt{43}) \approx 1.13476$. }
  \label{fig:0}
\end{figure}

 We have plotted the coefficients $\{ W_i, Z_i \}_{i=0}^4$ for $r \in (1, r^\ast )$ in Figure \ref{fig:0}. Note that at $r \approx r_n$ the singularity of $Z_n$ will propagate to every $W_i, Z_i$ for $i > n$ since they depend on $Z_n$ via the recurrence \eqref{eq:recurrence}.

\section{Barrier arguments}

\begin{figure}
\centering
  \includegraphics[width=.45\linewidth]{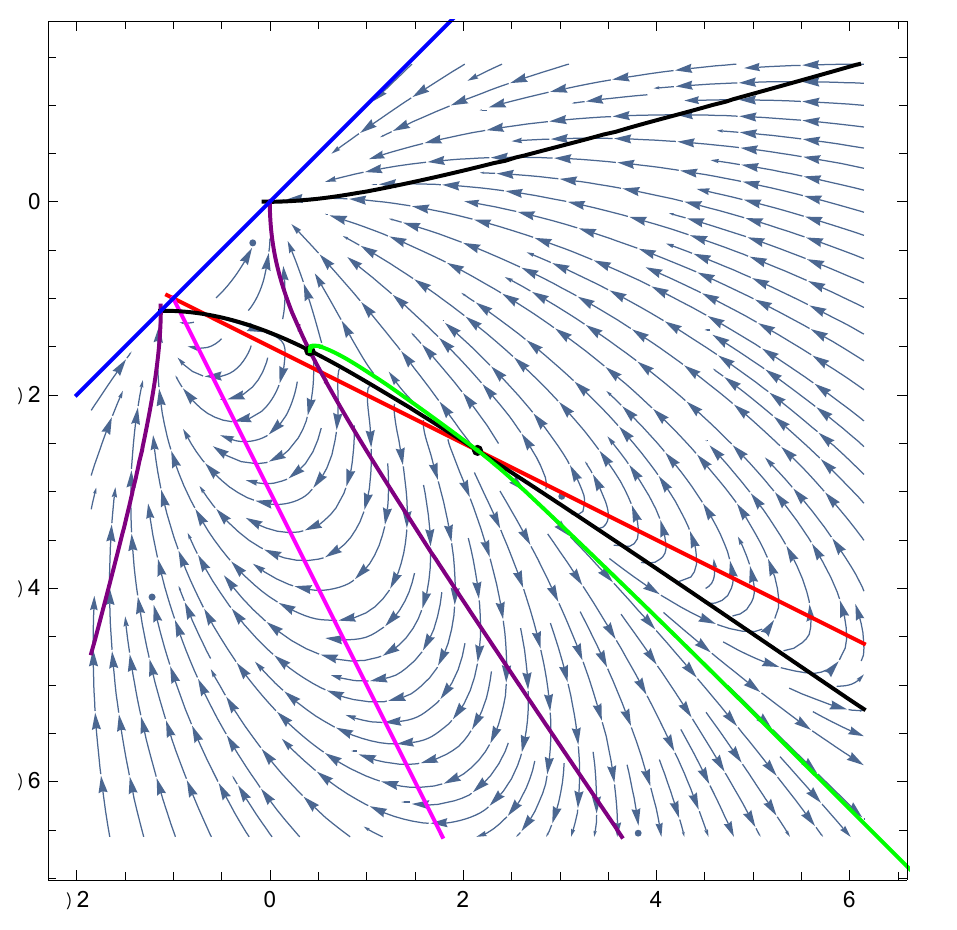}
  \includegraphics[width=.45\linewidth]{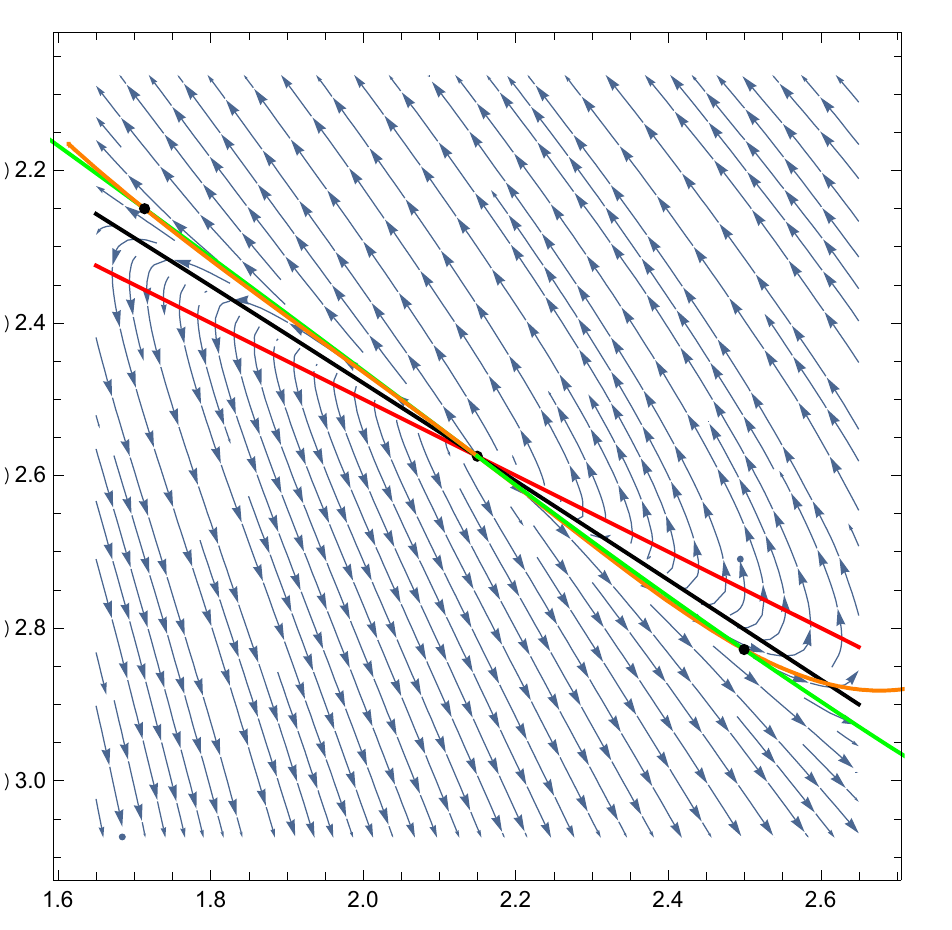}
  \captionof{figure}{\footnotesize Two plots at different scales of the vector field $(W_\psi, Z_\psi)$ for $r =1.13 \in (r_3, r_4)$. In the center we have $P_s$, which is the intersection of $D_Z = 0$ (red line) and $N_Z = 0$ (black curve). We also show $D_W = 0$ (in pink) and $N_W = 0$ (in purple). In green, we have the global barriers $b^{\rm{fl}}(t)$ (to the left of $P_s$) and $B^{\rm{fr}}(W, Z) = 0$ (to the right of $P_s$). At the smallest scale, we also show in orange the local barriers $b^{\rm{nl}}(t)$ (to the left of $P_s$) and $b^{\rm{nr}}_3(t)$ (to the right of $P_s$). For $b^{\rm{nr}}_3(t)$ we took $\beta = 500$. We have also indicated the points of intersection of the local barriers and the global barriers.  }
  \label{fig:1}
\end{figure}

For the sake of simplicity, let us concentrate on how to show items \ref{pt:1}, \ref{pt:2} and \ref{pt:3} for the case $\gamma = \frac53$ and $k \in (3, 4)$, which corresponds to $r \in (r_3, r_4)$. \\
% Right
The idea to prove item \ref{pt:1}\ is to construct two different barriers bounding the behavior of the solution. We will have one global barrier which we denote by $b^{\rm{fl}}(t)$ and a local one, which we denote by $b^{\rm{nl}}(t)$. They are given by
\begin{align*}
b^{\rm{nl}}(t)    &= \left( \sum_{i=0}^3 \frac{W_i}{i!} t^i, \sum_{i=0}^3 \frac{Z_i}{i!} t^i \right), \\
B^{\rm{fl}}(W, Z) &= \left( W_0 + B_1 W_1 t + \frac{B_2}{2} t^2, Z_0 + B_1 Z_1 t + \frac{B_3}{2} t^2 \right),
\end{align*}
where the coefficients $B_2, B_3$ are chosen so that $b^{\rm{fl}}(1) = P_\eye$ and the coefficient $B_1$ is chosen so that the barrier matches one further order with the ODE at $P_\eye$. The point $P_\eye$ is defined as the intersection of $N_W = N_Z = 0$ on the $W > Z$ half-plane. One can see both barriers in Figure \ref{fig:1}.

The global barrier will connect the point $P_s$ with $P_\eye$ in such a way that all trajectories of the ODE traverse $B^{\rm{fl}}$ upwards. Once we show that the smooth solution stays above $b^{\rm{fl}}(t)$ it is easy to conclude it has to converge to $P_0$, however, $b^{\rm{fl}}(t)$ will not be well-adapted to the geometry of the phase portrait at $P_s$, and this means that the smooth solution will not start above $b^{\rm{fl}}(t)$. In order to solve that, we use the local barrier $b^{\rm{nl}}(t)$ which matches the smooth solution up to third order, and thus is well-adapted to the geometry of the phase portrait at $P_s$. In particular, the smooth solution will start above $b^{\rm{nl}}(t)$ and trajectories will traverse $b^{\rm{nl}}(t)$ upwards for a short period of time $t \in [0, t_v]$. Thus, if we show that $b^{\rm{nl}}(t)$ and $B^{\rm{fl}}$ intersect at some time $t \in (0, t_v)$, we will be done, since the concatenation of the two barriers will correctly bound the behavior of the solution. This is done via a computer-assisted proof which involves a careful desingularization as $k \rightarrow 3^+$.

%Left
We now describe how to prove items \ref{pt:2} and \ref{pt:3}. Let us consider $n \in \{ 3, 4 \}$ and define the local barrier
\begin{equation} \label{eq:defb}
b^{\rm{nr}}_n(t) = \left( \sum_{i=0}^n \frac{W_i}{i!} (-t)^i,  \sum_{i=0}^n \frac{Z_i}{i!} (-t)^i + \frac{(-1)^n \beta |Z_n|}{(n+1)!} (-t)^{n+1} \right), \qquad \text{ for } 0 \leq t \leq \beta |k-n|^{\frac{1}{n-1}}
\end{equation}
which matches until $n$-th order with the smooth solution at $P_s$. We then define $$P^{\rm{nr}}_n (t) =  \left( \p_\psi W (b^{\rm{nr}}_n(t)), \p_\psi Z (b^{\rm{nr}}_n(t)) \right) \wedge  b^{\rm{nr} \; \prime}_n (t) ,$$ so that the sign of $P^{\rm{nr}}_n (t)$ informs us if the solutions of the ODE are traversing $b^{\rm{nr}}_n(t)$ in the upwards direction (negative sign) or in the downwards direction (positive sign). A careful computation of $P^{\rm{nr}}_n(t)$ yields
\begin{equation} \label{eq:P}
P^{\rm{nr}}_n(t) = \frac{| N_{W, 0} D_{Z, 1} | \beta}{(n+1)!} |Z_n| (-t)^{n+1} + \frac{ |\p_z D_Z (P_s) N_{W, 0} | }{n! (n-1)!} Z_n^2 (-t)^{2n-1} + o_{r \rightarrow r_n} \left( \beta |Z_n| t^{n+1} + Z_n^2 t^{2n-1} \right)
\end{equation}

For item \ref{pt:2}, we set $n=3$. Comparing the $(n+1)$-th Taylor coefficients of \eqref{eq:defb} and the smooth solution, one can see that $b^{\rm{nr}}_{3}(t)$ will be above the smooth solution near $P_s$ provided $\beta$ is chosen sufficiently large. Moreover, from \eqref{eq:P}, the barrier $b^{\rm{nr}}_{3}(t)$ will bound the trajectory of the smooth solution up to $t \leq c \beta |k-n|$. We construct another barrier $B^{\rm{fr}}$ given in implicit form by the nullset of:
\begin{equation*}
B^{\rm{fr}}(W, Z) = (W-W_0 - \frac12 Z + \frac12 Z_0 )( W+Z - F_0 ) - F_1(W+Z-W_0-Z_0).
\end{equation*}
The values of $F_0$ and $F_1$ are chosen so that $B^{\rm{fr}}$ matches the subleading order terms of its Taylor expansion with the smooth solution both at $P_0$ and $P_s$. Concretely, $F_0 = -2(r-1)$ and $F_1 = (W_0 + Z_0 - F_0) \frac{Z_1/2 - W_1}{W_1+Z_1}$. We can define $P^{\rm{fr}}_3(t)$ in the same way as we defined $P^{\rm{nr}}_3(t)$ and we show with a computer-assisted proof that $P^{\rm{fr}}_3(t) > 0$. That is, solutions always traverse $B^{\rm{fr}}(W, Z) = 0$ downwards. We have plotted $b^{\rm{nr}}_3(t)$ and $B^{\rm{fr}}(W, Z) = 0$ in Figure \ref{fig:1}.

Finally, we want to show that the concatenation of $b^{\rm{nr}}_3(t)$ and $B^{\rm{fr}}_3$ yields a barrier that bounds adequately the global behavior of the smooth solution. To that end, it suffices to show that both barriers intersect at some time $t_i \in (0, c \beta (k-3))$, so that $b^{\rm{nr}}_3(t)$ remains valid up to $t_i$. The choice of $F_0, F_1$ guarantees that the Taylor expansions at $P_s$ of $b^{\rm{nr}}_3(t)$ and $B^{\rm{fr}}$ first differ at their second order coefficients. Comparing their Taylor series one can conclude that both barriers intersect at some $t_i \les |k-3|$, so $t_i \leq c\beta (k-3)$ taking $\beta$ sufficiently large.

\begin{figure}
\centering
  \includegraphics[width=.45\linewidth]{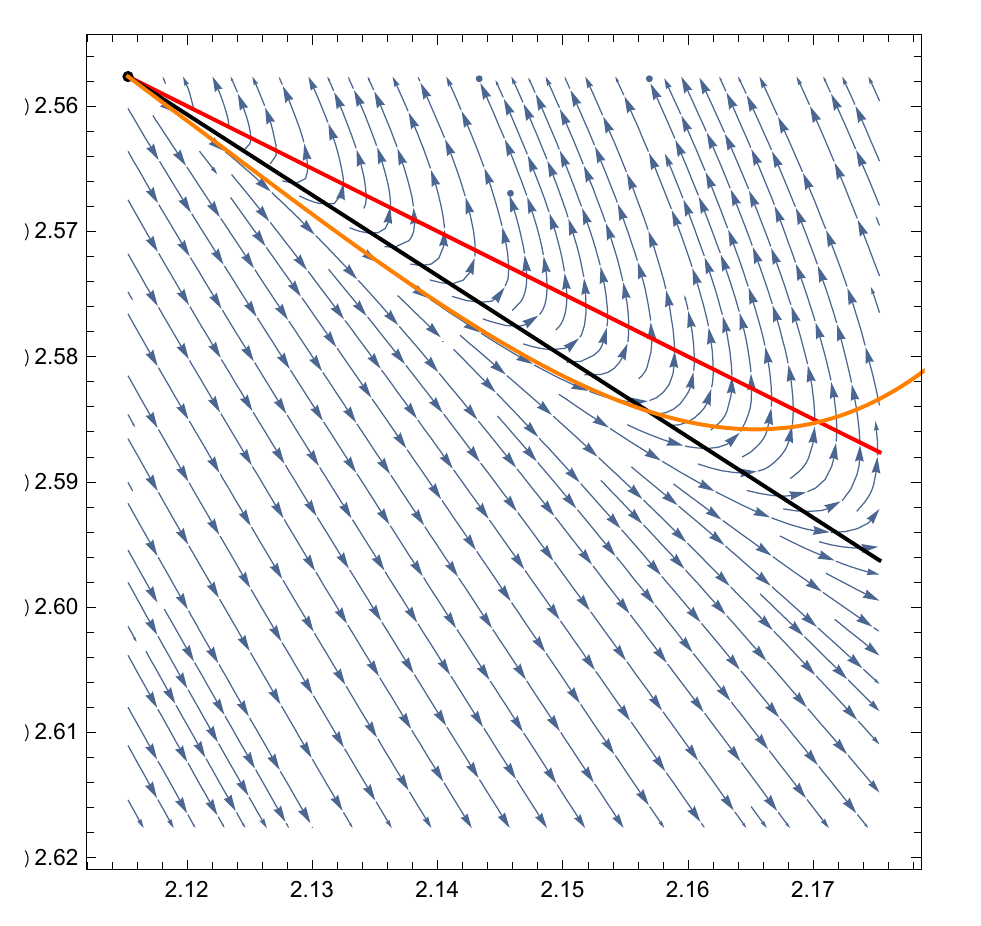}
  \captionof{figure}{\footnotesize A plot of the vector field $(W_\psi, Z_\psi)$ for $r =1.1347 \in (r_3, r_4)$ (moreover $1.1347 \approx r_4$). In the center we have $P_s$, which is the intersection of $D_Z = 0$ (red line) and $N_Z = 0$ (black curve). We also show $b^{\rm{nr}}_4(t)$ in orange, for $\beta = 500$. }
  \label{fig:2}
\end{figure}

With respect to item \ref{pt:3}, we set $n=4$ and use the local barrier \eqref{eq:defb}. Taking $\beta$ sufficiently large, $b^{\rm{nr}}_4(t)$ will be below the smooth solution for $t$ sufficiently small. Moreover, in this case, both terms from \eqref{eq:P} are negative, giving that solutions to the ODE cross $b^{\rm{nr}}_4(t)$ upwards for every $0 \leq t \leq \beta (4-k)^{1/3}$. Therefore, if we show that $b^{\rm{nr}}_4(t)$ intersects $D_Z = 0$ for some $t < \beta (4-k)^{1/3}$ we will be done. In that interval, we can compute
\begin{equation*}
D_Z (b^{\rm{nr}}_4(t)) = \frac{W_1 + 2Z_1}{3} t + \frac{Z_4}{36} t^4 + o_{k \rightarrow 4^-} ( (4-k)^{2/3} ),
\end{equation*}
where the two main terms are both of order $(4-k)^{1/3}$. Checking that the two main terms have different signs, we deduce that $b^{\rm{nr}}_4(t)$ intersects $D_Z = 0$ at some $t \leq \beta (4-k)^{1/3}$, provided that $\beta$ is chosen sufficiently large. We have a plot of this situation in Figure \ref{fig:2}.

\section{Computer-assisted proofs}

The representation of real numbers using a finite number of zeros and ones has the advantage of allowing finite calculations and a practical framework. However, this method also has the disadvantage of being limited to a finite (although large) amount of numbers and the potential for inaccuracies when performing mathematical operations. As an alternative, we will use upper and lower bounds for all relevant quantities, and propagate these bounds by rounding up or down as necessary to account for errors introduced by the computer during the calculation process.

We can now construct an arithmetic by the theoretic-set definition
\begin{align*}
[x] \star [y] = \{x \star y | \quad x \in [x], y \in [y]\},
\end{align*}
for any operation $\star \in \{+,-,\times,/ \}$. These are defined by the following equations:
\begin{gather*}
[x]+[y]  = [\nabla(\underline{x}+\underline{y}),\Delta(\overline{x}+\overline{y})], \quad 
[x]-[y] = [\nabla(\underline{x}-\overline{y}),\Delta(\overline{x}-\underline{y})],
\\
[x]\times[y]   = [\nabla(\min\{\underline{x}\underline{y},\underline{x}\overline{y},\overline{x}\underline{y},\overline{x}\overline{y}\}),
\Delta(\max\{\underline{x}\underline{y},\underline{x}\overline{y},\overline{x}\underline{y},\overline{x}\overline{y}\})], \\
[x] / [y]  = [x] \times \left[\frac{1}{\overline{y}},\frac{1}{\underline{y}}\right], \text{ whenever } 0 \not \in [y],
\end{gather*}
where $\nabla$ and $\Delta$ are respectively the round-down and round-up operators.

The main feature of the arithmetic is that if $x \in [x], y \in [y]$, then necessarily $x \star y \in [x] \star [y]$ for any operator $\star$. This property is fundamental in order to ensure that the true result is always contained in the interval we get from the computer. This process is completely rigorous and independent of the architecture or the software of the computer. We can also define functions of intervals $f([x])$. For example, if $f([x]) = [x] \times [x] + [x]$, then $f([-1,2]) = [-1,2] \times [-1,2] + [-1,2] = [-2,4]+[-1,2] = [-3,6]$.
%
%\color{blue} I don't think this is relevant in this paper. We don't integrate. \color{black}
%We now explain how to perform rigorous integration. Recall that from the Mean Value Theorem, it follows that if $f \in C^{1}$:
%\begin{align*}
% \int_{a}^{b}f(z)dz = f'(\xi)(b-a), \quad \xi \in [a,b] \Rightarrow 
% \int_{a}^{b}f(z)dz \in f'([a,b])(b-a).
%\end{align*}
%Following this principle and assuming enough regularity on $f$, one might develop more complicated integration schemes, such as for example
%\begin{align*}
%\int_{a}^{b} f(z) dz \in  \tfrac{b-a}{2}\left(f\left(\tfrac{b-a}{2}\tfrac{\sqrt{3}}{3} + \frac{b+a}{2}\right)+f\left(-\tfrac{b-a}{2}\tfrac{\sqrt{3}}{3} + \tfrac{b+a}{2}\right)\right)
%+\tfrac{1}{4320}(b-a)^{5}f^{4}([a,b]).
%\end{align*}
%
%Of course, all operations concerning the evaluation of $f$ and its derivatives are interval based, as described above, in order to guarantee the rigorous enclosure of the result. 

Early computer-assisted proofs were constrained to finite dimensional problems 
\cite{Fefferman-DeLaLLave:relativistic-stability-of-matter-I, 
Tucker:lorenz-focm}; however, recent advances in computational power  have 
enabled the methods to be adapted to infinite dimensional problems  (PDE). In the 
context of fluid mechanics we highlight the following equations: De 
Gregorio \cite{Chen-Hou-Huang:blowup-degregorio}, SQG 
\cite{Castro-Cordoba-GomezSerrano:global-smooth-solutions-sqg}, 
Whitham \cite{Enciso-GomezSerrano-Vergara:convexity-cusped-whitham}, Muskat \cite{GomezSerrano-GraneroBelinchon:turning-muskat-computer-assisted,Cordoba-GomezSerrano-Zlatos:stability-shifting-muskat-II},
Kuramoto-Shivasinsky 
\cite{Arioli-Koch:cap-stationary-ks,Figueras-DeLaLLave:cap-periodic-orbits-kuramoto,Gameiro-Lessard:periodic-orbits-ks,Figueras-Gameiro-Lessard-DeLaLLave:framework-cap-invariant-objects,Zgl,ZM1}, Navier-Stokes \cite{vandenBerg-Breden-Lessard-vanVeen:periodic-orbits-ns,Arioli-Gazzola-Koch:uniqueness-bifurcation-ns}, Burgers-Hilbert \cite{Dahne-GomezSerrano:highest-wave-burgershilbert} or the Hou-Luo model \cite{Chen-Hou-Huang:blowup-hou-luo}. We also refer the reader to the books \cite{Moore-Bierbaum:methods-applications-interval-analysis,Tucker:validated-numerics-book} and to the survey \cite{GomezSerrano:survey-cap-in-pde} and the book \cite{Nakao-Plum-Watanabe:cap-for-pde-book} for a more specific treatment of computer-assisted proofs in 
PDE.

In the paper \cite{implosion}, interval arithmetic is used to check the validity (positivity conditions) of the barriers and to compute a few thousands of coefficients of the Taylor expansion at $P_s$ (the latter is only used for the case $\gamma = \frac75$). We performed the rigorous computations using the Arb library
\cite{Johansson:Arb} and specifically its C implementation. 
The positivity checks involve using a branch and bound algorithm to evaluate the open conditions mentioned in the paper. We start by enclosing the condition within a box in a parameter space (which is at most 2-dimensional). If the enclosure provides a definite sign, we accept or reject it based on whether the sign matches the desired result. If the enclosure does not provide a sign, we split the box in half along one of the dimensions and repeat the process. This procedure continues until the maximum length in any dimension of the box reaches a tolerance of $10^{-10}$, at which point the program will fail. In our case, this tolerance was never reached. 

\section*{Acknowledgements}
T.B.\ was supported by the NSF grants DMS-2243205 and DMS-1900149, a Simons Foundation Mathematical and Physical
Sciences Collaborative Grant and a grant from the Institute for Advanced Study.  G.C.-L.\ was supported by a grant from the Centre de Formaci\'o Interdisciplin\`aria Superior, a MOBINT-MIF grant from the Generalitat de Catalunya and a Praecis Presidential Fellowship from the Massachusetts Institute of Technology. G.C.-L.\ would also like to thank the Department of Mathematics at Princeton University for partially supporting him during his stay at Princeton and for their warm hospitality.
This project has received funding from the European Research Council (ERC) under the European Union's Horizon 2020 research and innovation program through the grant agreement 852741 (G.C.-L.\, J.G.-S.). J.G.-S.\ was partially supported by NSF through Grant DMS-1763356 and by the AGAUR project 2021-SGR-0087 (Catalunya). J.G.-S.\ and G.C.-L.\ were partially supported by MICINN (Spain) research grant number PID2021--125021NA--I00.

\bibliography{euler}
\bibliographystyle{plain}

\end{document}